\documentclass[12pt]{amsart}
\usepackage{amscd,enumerate,amsfonts,calc,amsmath,verbatim}
\begin{document}
\title{Rings which are almost Gorenstein}
\author{Craig Huneke and Adela Vraciu}
\address {Department of Mathematics, University of Kansas, Lawrence,
KS 66045}
\email{huneke@math.ukans.edu}
\address{Department of Mathematics, University of South Carolina,
Columbia, SC}
\email{vraciu@math.sc.edu}
\date{\today}
\thanks{Both authors were partially supported by NSF grant
DMS-0098654}
\begin{abstract} We introduce classes of rings which are close to being Gorenstein.
These  rings arise naturally as specializations of rings of countable CM type. We
study these rings in detail, and along the way generalize an old result of Teter
which characterized Artinian rings which are Gorenstein rings modulo their socle.
\end{abstract}

\maketitle
\newcommand{\Hom}{\operatorname{Hom}_R}
\newcommand{\Ann}{\operatorname{Ann}}
\newcommand{\gr }{\operatorname{gr}}

\swapnumbers
\theoremstyle{plain}
\newtheorem{theorem}{Theorem}[section]

\newtheorem{prop}[theorem]{Proposition}
\newtheorem{lemma}[theorem]{Lemma}
\newtheorem{corollar}[theorem]{Corollary}
\newtheorem*{Corollary}{Corollary}

\theoremstyle{definition}
\newtheorem{note}[theorem]{Note}
\newtheorem{obs}[theorem]{Observation}
\newtheorem{definition}[theorem]{Definition}
\newtheorem*{Definition}{Definition}
\newtheorem{example}[theorem]{Example}
\newtheorem*{notation}{Notation}
\newtheorem*{conj}{Conjecture}
\newtheorem*{claim}{Claim}
\newtheorem*{question}{Question}

\newcommand{\li}{\tilde}
\newcommand{\aaa}{\mathfrak{a}}
\newcommand{\bbb}{\mathfrak{b}}
\newcommand{\ccc}{\mathfrak{c}}
\newcommand{\ub}{\underline{b}}

\newcommand{\m}{\mathfrak{m}}
\newcommand{\param}{\underline{x}}
\newcommand{\tpar}{\underline{x}^{[t]}}
\newcommand{\tparq}{\underline{x}^{[tq]}}
\newcommand{\bs}{\boldsymbol}
\newcommand{\tx}{\noindent \textbf}
\newcommand{\ld}{\ldots}
\newcommand{\cd}{\cdots}
\newcommand{\q}{^{[q]}}
\newcommand{\cor}{^{<q>}}
\newcommand{\spec}{I^{*sp}}
\newcommand{\f}{(f_1, \ld, \hat{f_i}, \ld, f_n)}
\newcommand{\eqq}{\Leftrightarrow}
\newcommand{\fs}{(f', f_1, \ld, \hat{f_i}, \ld , f_n)}
\newcommand{\ins}{I_1\cap \ld \cap \hat{I_i} \cap \ld \cap I_n}
\newcommand{\xs}{x_1, \ld, x_d}
\newcommand{\xt}{(x_1^t, \ld, x_d^t)}
\newcommand{\xtq}{(x_1^{tq}, \ld, x_d^{tq})}
\newcommand{\product}{x_1\cdots x_d}

\newcommand{\ui}{\underline{i}}
\newcommand{\und}{\underline}
\newcommand{\arr}{\Rightarrow}
\newcommand{\lar}{\longrightarrow}
\newcommand{\inc}{\subseteq}
\newcommand{\w}{\omega}

\newcommand{\Tor}{\mathrm{Tor}^R}
\newcommand{\Ext}{\mathrm{Ext}_R}
\newcommand{\x}{\underline{x}}
\newcommand{\ol}{\overline}
\newcommand{\syz}{\mathrm{syz}}
\newcommand{\Soc}{\mathrm{Soc}}

\section*{Introduction}
\bigskip

This paper began with a desire to understand better Cohen-Macaulay rings of countable or finite
representation type. Let $(R, \m)$ be a (commutative Noetherian) local ring of dimension
$d$.  Recall that a nonzero $R$-module $M$ is called {\it maximal
Cohen--Macaulay} (MCM) provided it is finitely generated and there
exists an $M$-regular sequence $\{x_1, \ldots, x_d\}$ in the maximal
ideal $\m$.

\begin{Definition} A Cohen-Macaulay local ring $(R, \m)$ is said to have {\rm
finite ({\it resp.},
countable) Cohen--Macaulay type} if it has only finitely (resp.,
countably) many isomorphism classes of indecomposable maximal Cohen--Macaulay modules.
\end{Definition}

A particular question we were interested in answering was the following:
what are the possible Hilbert functions of  $R/I$, where $R$ is a Cohen-Macaulay
ring of at most countable CM type, and $I$ is generated by a general system of parameters?
While we have not answered this question, what we found instead was that such
quotients behave much like Gorenstein rings in a very precise sense, and this
changed the direction of our inquiry to understanding these `almost' Gorenstein
rings. Our first section presents the basic formulation of the properties of
such rings. Our third section is devoted to proving that Cohen-Macaulay
rings having at most countable CM type are almost Gorenstein in the sense of having some of these
properties. Let $R$ be a local Cohen-Macaulay ring with canonical module $\w$.
An example of one of the properties we are considering is that
$\w^*(\w)$, the set of all elements of the form $f(x)$ where $x\in \w$
 and $f:\w\lar R$ is an $R$-linear map, contains the maximal ideal of $R$. If $R$
is Gorenstein, then of course $\w$ is free and $\w^*(\w)$ is the whole ring.

After introducing these rings, it is natural to look for examples.
It turns out that Artinian Gorenstein rings modulo their socle always are
examples, and this led us to a 1974 result of Teter \cite{Te} which gave
an intrinsic characterization of such rings. We are able to improve his
result if $2$ is a unit, by removing a  seemingly important technical assumption
of Teter's. This work is in Section 2. We further improve the result in
the case $2$ is a unit in the graded case, giving necessary and sufficient
conditions for a standard graded Artinian ring to be the homomorphic image
of a standard graded Gorenstein ring by its socle element.

The properties which characterize the rings we are interested in are all
closely related and perhaps are all equivalent: we have been unable to
decide whether they are equivalent, except in special cases. In section four, we prove that all
the conditions introduced in the first section are equivalent for Artinian
local rings of type two.

Finally our last section classifies $\m$-primary monomial ideals of type
three which are `almost' Gorenstein. An analysis of this classification shows that the conditions under consideration are equivalent in this case.

\section{Almost Gorenstein Rings}
\bigskip

In this section we introduce rings which are almost Gorenstein in a sense which the
first proposition identifies. We begin by
identifying several Gorenstein-like properties which a ring may have.
If $M$ is a module over a ring $R$, we denote $\Hom(M,R)$ as $M^*$, and
by $M^*(M)$ we mean the ideal consisting of all $f(x)$ where $x\in M$ and
$f\in M^*$.

\begin{prop}\label{almgor}
Let $(R,\m,k)$ be a Noetherian Cohen-Macaulay local ring with infinite residue field
$k$, and having a canonical module $\w$. Consider the following conditions:
\begin{enumerate}[\quad\rm (1)]
\item $\m\inc \w^*(\w)$.
\item For all ideals $K$ generated by a system of parameters and for all ideals
$I\supseteq K$, $K:(K:I)\inc I:\m$.
\item For all $K$ generated by a system of parameters which are not contained in a given
finite (or countable if $R$ is complete or $k$ is uncountable) set of primes not
equal to the maximal ideal, and for all ideals
$I\supseteq K$, $K:(K:I)\inc I:\m$.
\end{enumerate}

Then (1) $\implies$ (2) $\implies$ (3).
\end{prop}

{\bf Note:} In case $R$ is Artinian, the assumption on $K$ in condition (3) is void, thus in this case conditions (2) and (3) are automatically the same.

\begin{proof} Obviously (2) $\implies$ (3), so we need only to prove that (1) $\implies$ (2).

It is well-known that $a \in \omega ^*(\omega)$ if and only if the map $m_a:\omega \rightarrow \omega
$ given by multiplication by $a$ factors through a free module $F$ (\cite{D}, Lemma 1.2, Lemma 1.3).
 Recall that by Matlis duality we have $0:_R(0:_{\omega}I)=I$ for any ideal $I$. 
Consider the image of the submodule $N=0:_{\omega}I \subset \omega$. On
one hand, it is $a (0:_{\omega} I)$. On the other hand, the image in
$F$ is contained in $(0:_RI)F$, and the image of this latter module in
$\omega$ is contained in $(0:_RI)\omega$.

Therefore we have 
$$a(0:_{\omega}I)\subset (0:_R I)\omega;$$ taking duals in $R$, we get
$$
0:_Ra(0:_{\omega}I)\supset 0:_R(0:_R I)\omega
$$
But $0:_Ra(0:_{\omega}I)=I:a$, and $0:_R(0:_R I)\omega=0:_R(0:_R I)$
since $\omega $ is faithful.
\end{proof}

A basic question we are unable to answer except in some cases is:

\begin{question}
If $(R, \m)$ is an Artinian local ring which satisfies (3) of Proposition~\ref{almgor},
is it true that $\omega ^*(\omega )=\m$?
\end{question}

In other words, are the three properties of the above proposition equivalent?
\begin{obs}
We think of all of these properties as describing rings which are
almost Gorenstein. Of course, if $R$ is Gorenstein, then $K:(K:I) = I$, and
$\w^*(\w) = R$.  We note that Ding \cite{D} studied rings which satisfy property (1),
and we shall use some of his ideas in this paper.
\end{obs}
\smallskip

We now consider in some detail what condition (1) of Proposition~\ref{almgor} means.
Let $S$ be a Gorenstein Artininan ring, $K=(f_1, \ldots , f_n)\inc S$ an
ideal. We want to study the ring $R=S/(0:K)$. There is no loss of generality in writing
$R$ in this way, since every ideal in a Gorenstein Artinian ring is an annihilator ideal.
Note that the canonical module for this ring can be identified with $\text{Hom}_S(R,S)\cong
0:_S(0:_SK) = K$, so that $\omega_R \cong K$.

To study the trace ideal $\omega^*(\omega)$, consider an $R$-linear
map
$\phi:K \rightarrow S/(0:K)$.
Let $u_i = \phi (f_i)$, and let $v_i$ denote a lifting of $u_i$ to
$S$. We must have $u_i (0:f_i)=0$ in $R$, thus $v_i (0:f_i)\inc
0:K$ (in $S$). Since $S$ is Gorenstein, this means that $v_i \in
(0:K):(0:f_i) = (f_i):K$.
Hence the image of any $R$-linear map $\phi $ is contained in
$$
\frac{(f_1):_SK + \cdots + (f_n):_SK}{0:K}
$$
Thus, a necessary condition for $\omega ^*(\omega)=\m$ is that
\begin{equation}\label{star} (f_1):_SK + \cdots + (f_n):_SK=\m   \end{equation}
for any choice of a system of
generators $f_1, \ldots, f_n$ of $K$.

In particular, $\omega ^*(\omega )=\m$ implies that
\begin{equation}\label{2star}
\mathrm{there\  exists}\  i=1, \ldots, n\  \mathrm{such\ that\ }  (f_i):K \not\inc \m^ 2.
\end{equation}

The next result shows that this last property also holds under the
weaker assumption that $R$ satisfies property (3) of Proposition~\ref{almgor}.
\smallskip

\begin{prop}
Let $S$ be a Gorenstein local Artinian ring.
If $R=S/(0:K)$ is an Artinian local ring which satisfies property (3) of Proposition~\ref{almgor},
and if $f_1, \ldots, f_n$ is a minimal set of
generators for $K$, then there exists an $i=1, \ldots, n$ such that
$(f_i):K \not\inc \m^2$.
\end{prop}

\begin{proof}
Assume that $R$ satisfies property (3) of Proposition~\ref{almgor}, but $f_i :K \inc \m^2$ for
all $i$. We may assume that $\m K\ne 0$. For if not, then $K$ is either $0$ or
has exactly one generators, a representative of the socle of $S$. In either case, the
result is trivial.
Property (3) states that for every ideal $I\supset
0:K$, we have $$(0:K):((0:K):I)\inc I:\m$$ (all colons are taken in
$S$), or equivalently $IK:K\inc I:\m $.
The equivalence follows because
$$(0:K):((0:K):I)=(0:((0:K):I):K=(I (0:(0:K)):K = IK:K.
$$

Take $I_i=0:f_i$. Since $f_i \in K$, $I_i \supset 0:K$.
We have
$K(0:f_i):K\inc 0:\m f_i$
for all $i=1, \ldots, n$, and therefore
$$
\bigcap\,  (K(0:f_i)) :K \inc \bigcap\,  (0:\m f_i)= 0:\m K
$$
or equivalently (by duality)
$$
\m K\inc 0:\left(\bigcap\,  (K(0:f_i)) :K\right) =K\left(0:\bigcap\,  (K(0:f_i))\right)=
$$
$$
K\left(\sum\,  (0:K(0:f_i))\right)= K\left(\sum \, (f_i):K\right).
$$
This contradicts the assumption that $(f_i):K\inc \m^2$ for all $i$.

\end{proof}

\section{Teter's rings}

In \cite{Te} in 1974, Teter characterized Artinian local rings $R$ which are of the form
$S/(\delta)$, where $S$ is a Gorenstein local Artinian ring, and $\delta$
generates its socle. We shall prove such rings satisfy (1) of Proposition~\ref{almgor}.
Teter's main theorem gives necessary and sufficient conditions for $R$ to be of this
form. His theorem states:

\begin{theorem} Suppose that $(R,\m,k)$ is local Artinian and let $E$ be an injective hull
of $k$. Then $R$ is a factor of a local Artinian Gorenstein ring by its socle if and only
if there exists an isomorphism $\phi: \m\lar \m^{\vee}$ satisfying $\phi(x)(y) = \phi(y)(x)$ for
all $x,y\in \m$, where $(\quad)^{\vee}$ denotes $\Hom(\quad,E)$.
\end{theorem}

An immediate corollary of this theorem is the following:

\begin{corollar} Suppose that $(R,\m,k)$ is local Artinian and is the factor of a local
Artinian Gorenstein ring by its socle. Then $R$ satisfies all of the conditions of
Proposition~\ref{almgor}.
\end{corollar}

\begin{proof} Taking the Matlis duals of the injection of $\m$ into $R$ gives a surjective
map from $E$ onto $\m^{\vee}$. Composing this surjection with the inverse of the isomorphism
$\phi$ in Teter's theorem gives a homomorphism $f:E\lar \m$ which is onto. As $E$ is
isomorphic to the canonical module of $R$, this proves that $\w^*(\w)$ contains $\m$.
\end{proof}

Our purpose in this section is to show that the condition that $\phi$ satisfies
$\phi(x)(y) = \phi(y)(x)$ for
all $x,y\in \m$ is basically unnecessary if $2$ is a unit. 
To do so we first need to prove some
preliminary remarks concerning an involution on $\w^*$.

For every $f \in \omega ^*$, we can define another linear map
$\tilde{f} \in \omega ^*$ as described below:

Let $x \in \omega$; consider the
map $\phi_{f, x} : \omega \rightarrow \omega $ defined by $\phi _{f,x}(y)=f(y) x$.
Since $\Hom (\omega , \omega )\cong R$, there exists
a unique $r_{f, x}\in R$ such that $\phi _{f, x}$ is multiplication by
$r_{f, x}$, i.e.
$$
f(y)x=r_{f, x} y
$$ for all $x, y \in \omega, f \in \omega ^*$.

\begin{Definition}\label{tilde}
Define $\tilde{f}:\omega \rightarrow R$ by $\tilde{f}(x)=r_{f,x}$.
\end{Definition}

It is not hard to check that $\tilde {f}$ is a linear map, and
moreover the mapping $\Phi: \omega ^* \rightarrow \omega ^*$ is linear.

The basic property which $\tilde{f}$ has is that for all $x,y\in \w$,
\begin{equation}\label{tildef} f(x)y = \tilde{f}(y)x \tag{*}\end{equation}
which follows because by definition $\tilde{f}(y)x  = r_{f,y}x = f(x)y$.

We summarize some of the properties of $\tilde{f}$:

\begin{prop} Let $(R,\m,k)$ be a local Artinian ring. Define $\Phi: \omega ^* \rightarrow \omega ^*$
as above. For $f\in \w^*$, set $I_f = f(\w)$, and
$J_f =\tilde{f}(\w)$.

\begin{enumerate}[\quad\rm (1)]
\item $\Phi $ is an isomorphism.
\item  $\Phi^2$ is the identity map on $\w^*$.
\item $\omega ^* (\omega)= \sum I_f = \sum J_f$.
\item $\ker f = 0:_{\omega} J_f$.
\item $\Hom (I_f,\omega ) \cong J_f $, i.e. $J_f$ is Matlis dual to $I_f$.
\end{enumerate}
\end{prop}

\begin{proof}
We claim that $\Phi $ is injective. Indeed, if $\tilde{f}=\tilde{g}$,
then $f(y)x= r_{f, x}y =r_{g, x}y=g(y) x$ for all $x, y \in \omega$. Since
$\omega $ is faithful, this shows that $f=g$. It follows that $\Phi$ is an
isomorphism since $\w^*$ has finite length and any injective map is also
surjective.

Part (2) follows at once from the basic property (*) above, which identifies
$r_{\tilde{f},x}y = \tilde{f}(y)x = f(x)y$, and so $r_{\tilde{f},x} = f(x)$,
which means that $\Phi^2$ is the identity map on $\w^*$.

Part (3) is clear from the definition.

To prove (4), consider
$0:_{\omega} J_f  = \lbrace y \in \omega \,  : \, r_{f, x} y =0 \,\text{for all}\, x\in \omega \rbrace =
\lbrace y \in \omega \,  : \, f(y)x =0 \, \text{for all}\,  x \in \omega \rbrace = \ker(f)$
(the last equality is  because $\omega $ is faithful).

Thus, we have
$$
I_f = f(\w) \cong \frac{\omega}{\ker (f)} = \frac{\omega}{0:_{\omega} J_f}
$$

Finally, we prove (5):
$$\Hom\left(  \frac{\omega}{0:_{\omega} J_f}, \omega \right) \cong
\lbrace \phi : \omega \rightarrow \omega : \phi (0:_{\omega} J_f)=0
\rbrace \cong 0:_R (0:_{\omega} J_f) =J_f
$$
\end{proof}

The following observation illustrates another connection between the
first and third properties of our basic proposition.

\begin{obs}
$$
\frac{0:(0:I)}{I}
$$ is the kernel of the canonical map
$$
\frac{R}{I} \rightarrow \left(\frac{R}{I}\right)^{**}.
$$
Thus property (3) in Prop. ~\ref{almgor} is satisfied if and only if the kernel of this canonical map is contained in the socle of $R/I$ for every ideal $I$.

On the other hand, $\omega^*(\omega)=\m $ if and only if the kernel of
the canonical map $\omega \rightarrow (\omega)^{**}$ is contained in the socle of
$\omega $.
\end{obs}
\begin{proof}
To see the first claim, note that $(R/I)^* = 0:I$, and consider the short exact sequence
$$
0\lar 0:I \lar R \lar \frac{R}{0:I}\lar 0
$$
Applying $\Hom(\quad , R)$ yields an exact sequence
$$
0 \lar \Hom \left( \frac{R}{0:I},\,  R\right) \lar R \lar \Hom(
0:I,\,  R)
$$
Since the module on the left is $0:(0:I)$, and the module on the right
is $(R/I)^{**}$, we obtain an exact sequence
$$
0 \lar 0:(0:I) \lar R \lar \left(\frac{R}{I}\right) ^{**}
$$
which proves the claim.

To see the second claim, notice that the kernel of the map $\omega
\rightarrow (\omega )^{**} $ is
$$ \bigcap_{f\in \omega ^*}
\ker(f)
= \bigcap 0:_{\omega } I_{\tilde{f}} = 0:_\omega  \left( \sum _{f \in \omega
^*}  I_f \right) = 0:_{\omega } \omega^*(\omega )
$$
\end{proof}

We now give our improvement of Teter's theorem. The use of the map $\tilde{f}$
makes it possible to avoid the awkward hypothesis in Teter's theorem that the
isomorphism $\phi$ satisfy $\phi(x)(y) = \phi(y)(x)$.

\begin{theorem}\label{genteter}
Let $(R, \m, k)$ be an Artinian ring with canonical module
$\omega$. Assume that $2$ is a unit in $R$ and $\Soc(R)\inc \m^2$.
Then $R$ is the quotient of a zero dimensional Gorenstein ring by its
socle if and only if there exists a surjective map
$f:\omega \rightarrow \m$.
\end{theorem}

\begin{proof} First assume that $R$ is the quotient of a zero-dimensional Gorenstein ring
by its socle. The result then follows from Teter's theorem, but the proof is so direct that
we give it here for the reader's convenience.
Let $R = S/(\delta)$, where $S$ is Gorenstein and $\delta$ generates the socle of $S$.
The exact sequence $$0\lar k\lar S\lar R\lar 0$$ gives, upon dualizing into the injective
hull of the residue field of $S$, the exact sequence,
$$0\lar E\lar S\lar k\lar 0$$
which proves that the injective hull $E$ of the residue field of $R$ is isomorphic to the
maximal ideal of $S$ (which is an $R$-module), and this clearly maps surjectively onto
the maximal ideal of $R$.

To prove the harder direction, we need only prove that given
$f$ as in the statement of the theorem, there exists an isomorphism $g:\m\lar  \m^{\vee}$ such that
$g(x)(y) = g(y)(x)$ for all $x,y\in \m$.

Set $h = (f+\tilde{f})$, which is a homomorphism from $\w\lar R$. We claim that
the kernel of this map is the socle of $\w$.  Suppose that $h(x) = 0$. Then
$(f+\tilde{f})(x) = 0$. By the definition of $\tilde{f}$ (see \ref{tildef}), it follows that
for all $y\in \omega$, $f(y)x = \tilde{f}(x)y = -f(x)y$, and hence
$$xf(y) + yf(x) = 0.$$

Since $f$ maps $\w$ onto $\m$, the kernel of $f$ must be the socle of
$\w$ (which is one-dimensional), and since $f(f(x)y-f(y)x)=f(x)f(y)-f(y)f(x)=0$ for all $x,y\in \w$, it follows that $\m (f(x)y-f(y)x) = 0$.

Suppose that $x \in \m\omega$. Write $x = \sum r_ix_i$, $r_i\in \m$.
Then
$$xf(y) + yf(x) = \sum r_i(x_if(y)+yf(x_i)),$$
but as $r_i\in \m$,
this sum is just $\sum r_ix_if(y) = 2xf(y)$. Hence $xf(y) = 0$.
Since $f(\w) = \m$, as we vary $y$ it follows that $x$ is in the socle of $\omega$. In order to finish the proof of the claim we need to show that $\ker (h)\inc \m \omega$. Let $x \in \ker (h)$. It follows that $\m x \inc \mathrm{ker } (h)\cap \m \w \inc \mathrm{Soc}(\w)$, and thus $x \in 0:_{\w} \m^2$. The assumption that $\mathrm{Soc}(R)\inc \m^2$ is equivalent, by Matlis duality, to $0:_{\w} \m^2\inc \m \w$.

Since the kernel of $h$ is $1$-dimensional, the length of the image of $h$ is exactly
one less than the length of $R$. It follows that $h$ maps $\w$ onto $\m$.

Next we prove that if $x,y\in \w$, then $h(x)y = h(y)x$. For
$h(x)y = (f+\tilde{f})(x)y = (f(x)y + \tilde{f}(x)y) =
(f(x)y + f(y)x)$. Since this is symmetric with respect to
$x$ and $y$, the claim follows.

Taking Matlis duals, we get a map $h^{\vee}: \m^{\vee}\lar R$ which is injective,
and so $h^{\vee}$ is an isomorphism of $\m^{\vee}$ and $\m$. Set $g$ equal to the
inverse map of $h^{\vee}$. Then $g$ is an isomorphism of $\m$ and $\m^{\vee}$. We
claim that for all $x,y\in \m$, $g(x)(y) = g(y)(x)$. We can then apply Teter's theorem
to finish the proof. The homomorphism $g$ is the dual of the inverse of $h$, where we
think of $h$ as an isomorphism of $\w/(\delta)$ with $\m$. If $u,v\in \m$, write
$u = h(x)$ and $v = h(y)$. Then $h^{-1}(u)v = xv = xh(y) = h(x)y = uh^{-1}(v)$, proving
that $h^{-1}$ satisfies the same symmetry condition. Taking Matlis duals preserves this
condition, proving the theorem.
\end{proof}

In the case $R$ is graded we can do better. In particular, the assumption that
the socle of $R$ be contained in $\m^2$ can be removed in the sense that if it is not
true, the structure of $R$ is fixed. Recall a Noetherian graded ring $R$ is
standard graded if $R_0 = k$ is a field and $R = k[R_1]$. 

\begin{theorem}\label{gengradteter}
Let $R$ be an Artinian standard graded ring over a field $k$, not having
characteristic $2$, with graded canonical module
$\omega$.  Set $\m$ equal to the ideal generated by all elements of positive degree.
The following are equivalent:

\smallskip
(1) Either $R\cong k[X_1,...,X_n]/(X_1,...,X_n)^2$, or $\Soc(R)\inc \m^2$
and there exists a degree $0$ graded surjective homomorphism $f:\omega(t) \rightarrow \m$ for some $t$.

\smallskip
(2) $R$ is the quotient of an Artinian standard graded Gorenstein ring by its
socle.
\end{theorem}

\proof First assume (2). If $\Soc(R)$ is not contained in $\m^2$, then there
must be a socle element of degree $1$, say $\ell$. Hence lifting back to $S$,
it follows that $\m_S\ell\inc \Soc(S)$, so that the socle of $S$ must live in
degree $2$. It follows that the Hilbert function of $S$ is $1, n, 1$ for some $n$,
and hence the Hilbert function of $R$ is $1, n$, implying that 
$R\cong k[X_1,...,X_n]/(X_1,...,X_n)^2$. Suppose that $\Soc(R)\inc \m^2$.
The graded canonical module of $S$ is $S(t)$, for some $t$, and then the
graded canonical module of $R$ is Hom$_S(R,S(t)) = \text{Hom}_S(R,S)(-t) = \m_S(-t)$.
Hence there is a graded surjective map onto $\m_R$ after twisting by $t$.

Conversely, we construct $S$ explicitly assuming (1). First, if
$R\cong k[X_1,...,X_n]/(X_1,...,X_n)^2$, then we may take
$$S = k[X_1,...,X_n]/(X_iX_j, X_i^2-X_j^2)$$ where the indices range over
all $1\leq i < j\leq n$. 

Otherwise, assume that $\Soc(R)\inc \m^2$, and there exists a degree $0$ graded
surjective homomorphism $f:\omega(t) \rightarrow \m$ for some $t$.
Define a ring structure on $S=k\oplus \omega(t)$ with multiplication
$$
(\alpha _1, x_1) (\alpha_2, x_2)=(\alpha _1 \alpha _2,
\alpha_1x_2+\alpha_2x_1 +\frac{x_1f(x_2) +x_2f(x_1)}{2})
$$
This multiplication is obviously commutative and gives $S$ a graded
structure. Moreover, $S$ is standard graded; the surjection of $\w (t)\rightarrow \m$ must have
a kernel of length $1$, by counting lengths. Hence the kernel is exactly the socle of
$\w$, and the socle is never a minimal generator of $\w$ unless $\w = R = k$.
Hence the minimal generators of $\w(t)$ correspond to the minimal generators of $\m$ and
all have degree $1$. In order to check
associativity, observe that the kernel of $f$ must be $\Soc(\omega)$,
since it is an $R$-submodule of $\omega$ of length one.
This implies that $uxf(y)=uyf(x)$ for any $u \in \m$, $x, y \in
\omega$.

$S$ is a graded ring with homogeneous maximal ideal $0 \oplus \omega$ (note that
$(0,x)^n =(0, x^nf(x))$, and since $x$ is nilpotent in $R$, $(0, x)$
is nilpotent in $S$).

Let us compute the socle of $S$. If $(\alpha, x )\in \Soc(S)$, we must
have $\alpha=0$. Thus,
$$\Soc(S)=\lbrace (0, x)\, |\,  x\in \omega , xf(y)+yf(x)=0\,  \forall y \in
\omega \rbrace = 0 \oplus \mathrm{Soc}(R)
$$
as seen in the proof of the theorem.

Note that $R$ is indeed a quotient of $S$ by its socle, since we have
a surjective map $S \rightarrow R$ given by $(\alpha , u )\rightarrow
\alpha + f(u)$. It is not hard to check that this is a ring
homomorphism; surjectivity is obvious, and by dimension counting it follows that
the kernel of the this map is the socle of $S$. \qed

The above proof works also in the non-graded case, providing a very different
approach than in the paper of Teter.

\bigskip
\section{Rings of Finite or Countable Cohen-Macaulay Type}
\medskip
Let $(R, \m)$ be a Cohen-Macaulay ring of countable Cohen-Macaulay
type, i.e. there are only countably many isomorphism classes of
indecomposable maximal Cohen-Macaulay modules (MCM's). Let $d$ denote
the dimension of $R$.

Let $\lbrace M_i \rbrace _i $ be a complete list of all of non-isomorphic indecomposable
MCM's (up to isomorphism); consider the set $\Lambda$ consisting of all the
annihilators of modules of the type:
$\mathrm{Ext} ^i_R(M_j, M_k)$, for $i\ge 1$.

Assume that the residue field is uncountable. Vector field
arguments show that $\m$ is not
contained in any countable union of proper subideals, in particular it is not contained in the union of all ideals in
$\Lambda $ other than $\m $ itself. Consider an element $x\in R$ not in
the union of the ideals in $\Lambda$ other than possibly $\m$. We will call the elements $x \in R$ satisfying this condition {\it general elements}. By a
{\it general system of parameters} $\x$ we mean a system of parameters such
that $(\x)$ contains a general element. If a general $x$ annihilates any
one of the modules listed above, it follows that $\m $ annihilates
that module.

\begin{prop}\label{ann}
The following modules have annihilator in $\Lambda$:
\begin{enumerate}[\quad\rm (1)]
\item $\Ext^i(M, N)$ for $i \ge 1$, $M$ MCM, $N$ arbitrary
\item $\Ext^i(N, M)$ for $i\ge d+1$, $M, N$ arbitrary
\item ${\mathrm Ext}_{\ol{R}}^i(N, M/(\x)M)$ for $i \ge 1$, $M$ MCM, $N$ arbitrary, and
$\x$ a s.o.p. such that $(\x)N=0$, where $\overline{R} = R/(\x)$.
\end{enumerate}
\end{prop}

\begin{proof}
(1) Let $N$ be an arbitrary module. Consider the Cohen-Macaulay
approximation of $N$ (see \cite{AB}):
$$
0 \lar C\lar T \lar N\lar 0
$$
where $C$ has finite injective dimension and $T$ is MCM.
Applying $\Hom (M, \quad )$ and using the fact that $\Ext ^i(M, C)=0$ we
get $\Ext^i(M, N)\cong \Ext^i(M, T)$ for $i\ge 1$.

(2) If $i \ge d+1$, we have $\Ext^i(N, M)\cong \Ext^{i-d}(\syz^d N,
M)$, and the result follows from (a) since $\syz^dN$ is MCM.

(3) We have Ext$^i_{\overline{R}}(N, M/(\x)M)\cong \Ext^{i+d}(N, M)$.

If $R$ is Gorenstein, the same proof in (a) works to show that
$\Tor_i(M, N)\cong \Tor_i(M, T)$, since in this case $C$ has finite
projective dimension.
\end{proof}

\begin{prop}\label{frac}
Let $M$ be a MCM, $(\x)\inc \aaa \inc I$ ideals in $R$, and $\x$
a general sop.
Then $\m$ annihilates
$$
\frac{ (\x)M:(\aaa :I)}{(\x)M+I(\x M:\aaa)}
$$
\end{prop}
\begin{proof}
Let $I=(f_1, \ldots, f_n)$. Consider the short exact sequence
$$
0 \lar \frac{R}{\aaa:I} \lar \bigoplus \frac{R}{\aaa}\lar N\lar 0
$$
where the first map is $\overline{u}\lar (\ol{f_1u}, \ldots, \ol{f_nu})$
Apply Hom$_{\overline{R}}( \quad , M/(\x)M)$. By Proposition~\ref{ann} (3) we know that $\m$
annihilates Ext$^1_{\overline{R}}(N, M/(\x)M)$, and therefore it annihilates the
cokernel of the induced map
$$
\bigoplus \text{Hom}_{\overline{R}}\left( \frac{R}{\aaa}, \frac{M}{\x M} \right) \lar
\text{Hom}_{\overline{R}}\left( \frac{R}{\aaa :I}, \frac{M}{\x M} \right)
$$
which is equivalent to
$$
\bigoplus \frac{\x M: \aaa}{ \x M}\lar \frac{\x M:(\aaa :I)}{\x M}.
$$
The above map is given by
$(\ol{u_1}, \ldots, \ol{u_n}) \lar \ol{f_1u_1+\ldots + f_nu_n}$, and
therefore the cokernel is
$$
\bigoplus \text{Hom}_{\overline{R}}\left( \frac{R}{\aaa}, \frac{M}{\x M} \right) \lar
\text{Hom}_{\overline{R}}\left( \frac{R}{\aaa :I}, \frac{M}{\x M} \right)
$$
which is equivalent to
$$
\bigoplus \frac{\x M: \aaa}{ \x M}\lar \frac{\x M:(\aaa :I)}{\x M}.
$$
The map is given by
$(\ol{u_1}, \ldots, \ol{u_n}) \lar \ol{f_1u_1+\ldots + f_nu_n}$, and
therefore the cokernel is
$$
\frac{\x M:(\aaa :I)}{\x M + I(\x M:\aaa)}.
$$
\end{proof}

\begin{corollar}\label{almostG}
If $\x $ is a general sop, we have $\m (\x :( \x :I))\inc (\x, I)$
for any ideal $I$.
\end{corollar}
\begin{proof}
Take $M=R$, $\aaa =(\x )$ in Proposition~\ref{frac}.
\end{proof}

Corollary~\ref{almostG} shows that Cohen-Macaulay local rings of finite or countable
CM type satisfy the third property of Proposition~\ref{almgor}.
We believe that they also satisfy the first condition that $\m\inc \w^*(\w)$, but have been
unable to prove it in this generality.

\bigskip
\section{Type Two Ideals}
\medskip
In this section we prove that the conditions of Proposition~\ref{almgor} are
equivalent for Artinian rings of type two. 

We begin with a general result concerning type $n$ ideals. The following notation will be used throughout the rest of the paper:

\begin{notation}\label{notation}
Let $(S, \m)$ be a local Gorenstein ring, and let $R=S/I$, where $I$ is an $\m-$primary ideal of type $n$.
We represent $I$ as an irredundant intersection of $n$ irreducible ideals $I= J_1 \cap J_2 \ \ldots \cap J_n$,
and choose $J\inc I$ be an irreducible ideal. Then for every $i=1, \ldots, n$ there exists $f_i \in S$
such that $J_i=J:f_i$.
\end{notation}

Since $I = J_1\cap ...\cap J_n$, then $J:I = \sum_i J:J_i$. To see this it suffices
to prove equality after computing annihilators into $J$. But 
$I = J:(J:I)$, while $J: (\sum_i J:J_i) = J:(J:J_1)\cap ...\cap J:(J:J_n) = J_1\cap ...\cap J_n = I.$
It follows that the canonical module of $S/I$ can be computed as:
\begin{equation}\label{canonical}
\omega _{S/I}\cong 
 \frac{J:I}{J} =
\frac{\sum_i (J:J_i)}{J}=
\frac{(J, f_1, \ldots, f_n)}{J}.
\end{equation}
\begin{prop}\label{genlcase} Adopt the notation above.
A necessary condition for $R = S/I$ to satisfy the third condition of Proposition~\ref{almgor} is that
for all $1\leq i\leq n$
$$
\sum_{j\ne i}(J_i:J_j) +  (J_1:J_i) \cap \ldots \cap (J_n :J_i) = \m
$$
\end{prop}

\begin{proof}

Without loss of generality we prove the claim for $i = 1$. Set $K = J_2\cap ... \cap J_n$ and $J_1 = J$.
The asserted equality is equivalent to saying that $K:J + J:K = \m$.
We may assume that $I = 0$. The third condition of  Proposition~\ref{almgor} tells us
that $0:(0:J)\inc J:\m$. We analyze $0:(0:J)$:
$$0:(0:J) = (K\cap J):((K\cap J):J) = (K\cap J): (K:J)$$
$$= K:(K:J)\cap J:(K:J).$$
The ideal  $K:(K:J)$ contains $J$, and so since $J$ is Gorenstein, it can be written
in the form $J:L$ for some ideal $J\inc L$. Moreover, since $K:(K:J)$ contains $K+J =
J:(J:K)$, we obtain that $L\inc J+K$. We then have that
$$0:(0:J) = J:L\cap J:(K:J) = J:(L+(K:J)).$$
The assumption that  $0:(0:J)\inc J:\m$ then implies that $J:(L+(K:J))\inc J:\m$ and
so $\m = J:(J:\m)\inc J:(J:(L+(K:J))) = L + K:J$. It follows that $J:K + K:J = \m$,
which gives the required formula. 
\end{proof}

\begin{theorem}\label{typetwo}
Let $(R,\m,k)$ be a local Artinian ring of type two. Write $R = S/I$,
where $S$ is a regular local ring, and write
$I=J_1 \cap J_2$, where $J_i$ are irreducible ideals.  Then the
three conditions of Proposition~\ref{almgor} are equivalent.
Furthermore, these conditions are also all equivalent to
$$J_1:J_2 +J_2:J_1 = \m.$$
\end{theorem}

\begin{proof} The fact that the weakest of the conditions in 
Proposition~\ref{almgor}, i.e. $0:(0:I)\inc I:\m$ for every ideal $I$ of $R$,
implies that
$J_1:J_2 +J_2:J_1 = \m$ is a particular case of Proposition ~\ref{genlcase}.

Assume $J_2:J_1 + J_1:J_2 = \m $ and we will prove that $\w^*(\w)$ contains $\m$.
Using the notation in (\ref{notation}),
 let  $y \in J_1:J_2 =(J, f_2):f_1$,
and choose $z$ such that
$yf_1 \equiv zf_2$ (mod $J$).

Recall the description of $\omega$ given in (\ref{canonical}), and define the  $R$-linear map $\Phi : \omega \rightarrow \frac{R}{I}$ by $\Phi (f_1)=z$,
$\Phi (f_2)=y$.
To check that $\Phi $ is well-defined, consider any relation $af_1 +
bf_2 \equiv 0$ (mod $J$). We need to check that $a z + by \in I
=J:(f_1, f_2)$.
Indeed, $af_1 z + bf_1 y = af_1z + b f_2z=(af_1+ bf_2) z\in J$, and
$af_2 z + bf_2 y = a f_1 y +bf_2 y =(af_1 + bf_2) y \in J$.
We've shown that $J_1:J_2 \inc \w^*(\w)$, and by symmetry the same holds
for $J_2:J_1$. It follows that $\m\inc \w^*(\w)$, completing the proof the theorem.

\end{proof}

Note that the above proof shows that all of these conditions are equivalent to
saying that $I:(I:J_1 )\inc J_1 :\m $. It is remarkable that this condition
implies the strongest condition, namely that $\w^*(\w)$ contains $\m$, especially
since it is not obviously symmetric in $J_1$ and $J_2$.

\begin{example}\label{type2mon} We can analyze type two monomial ideals completely. Let $S=k[x_1, \ldots, x_n]$ be a polynomial ring, and assume that $I$ is generated by monomials; it follows that $I$ can be represented as an intersection of irreducible monomial ideals.
It is well-known that the
only irreducible $\m$-primary monomials ideals are generated by powers of the variables.
Hence we may assume that
$$
J_1=(x_1^{c_1}, \ldots, x_n ^{c_n}), \ J_2 =(x_1^{d_1}, \ldots,
x_n^{d_n}), \
$$
$$
J_1 :J_2 = (x_1^{c_1}, \ldots, x_n ^{c_n}, x_1^{c_1-d_1}\cdots x_n^{c_n-d_n}), \text{and}$$
$$
J_2:J_1 =(x_1^{d_1}, \ldots, x_n^{d_n}, x_1^{d_1-c_1}\cdots x_n^{d_n-c_n})$$
where by convention if a variable has negative or zero exponent we drop it in the
product.
The only way that $J_1 :J_2 + J_2 :J_1 = \m$ is if one of the
following holds:

a. for every $i$
either $c_i=1$, or $d_i=1$.
By relabeling the variables, we can assume in this case that
$$J_1 = (x_1, \ldots, x_s, x_{s+1}^{c_{s+1}}, \ldots, x_n ^{c_n}),
J_2 =(x_1^{d_1}, \ldots, x_{s-1}^{d-1}, x_{s+1}, \ldots, x_n)$$
$$I=(x_1^{d_1}, \ldots, x_{s-1}^{d_{s-1}}, x_s^{c_s}, \ldots,
x_n^{c_n}, x_ix_j \, | \, i \le s, j \ge s+1)
$$

b. one of  $x_1^{c_1-d_1}\cdots x_n^{c_n-d_n}$ or $x_1^{d_1-c_1}\cdots x_n^{d_n-c_n}$
has all but one exponent non-positive, and the remaining
exponent (say the exponent of $x_n$) equal to one; moreover, for all
$i=1, \ldots, n-1$, either $c_i =1$ or $d_i =1$. Say that $c_1\le d_1,
\ldots, c_{n-1}\le d_{n-1}, c_n =d_n +1$. Then the second part of
the condition forces $c_1 = \cdots = c_{n-1} =1$.
In this case, the ideals are:
$$J_1= (x_1, \ldots, x_{n-1}, x_n^c),
J_2 = (x_1^{d_1}, \ldots, x_{n-1}^{d_{n-1}}, x_n^{c-1})
$$
$$
I=(x_1^{d_1}, \ldots, x_{n-1}^{d_{n-1}}, x_n^{c}, x_ix_n^{c-1}\, | \,
i=1, \ldots, n-1)
$$

c. both of $x_1^{c_1-d_1}\cdots x_n^{c_n-d_n}$ and  $x_1^{d_1-c_1}\cdots x_n^{d_n-c_n}$ have all but one exponent non-positive, and the remaining
exponent equal to one. For instance, $c_n-d_n=1$, $c_i\le d_i$ for all $i\ne n$,  and $d_{n-1}-c_{n-1}=1$, $d_i\le c_i$ for all $i\ne n-1$. This means that $x_n=x_1^{c_1-d_1}\cdots x_n^{c_n-d_n}, x_{n-1}=x_1^{d_1-c_1}\cdots x_n^{d_n-c_n}$. The conditions stated above imply that $c_i=d_i$ for all $i<n-1$; in order for $x_i$ to be in $J_1:J_2 +J_2:J_1$ we need $c_i=d_i=1$, thus $R$ is isomorphic to $S_0/I_0$, where $S_0=k[x, y]$ and $I_0=(x^c, y_d, x^{c-1}y^{d-1})$ for some $c, d$, thus $R$ is a Teter ring.
\end{example}

\bigskip
\section {Type Three Monomial Ideals}

\medskip

In this section we classify the type three primary monomial ideals $I$ in  a polynomial ring $S=k[x_1, \ldots, x_n]$ such that $R=S/I$ is Artinian and satisfies one  of the conditions of
Proposition~\ref{almgor}. We also prove that the three conditions are equivalent in this case.

The following notation will be assumed throughout this section:

$S=k[x_1, \ldots, x_n]$, $I=J_1\cap J_2\cap J_3$, where 
$J_i=(x_1^{a_{i1}}, \ldots, x_n^{a_{in}})$ is a monomial $\m-$ primary Gorenstein ideal (it is well-known that all monomial $\m-$ primary Gorenstein ideals are of this form).
We use $x^{a_i-a_j}$ to denote 
$\displaystyle \Pi x_k^{a_{ik}-a_{jk}}$, where the product runs over $k\in \{ 1, 2, \ldots, n\}$ with $a_{ik}\ge a_{jk}$. Note that we have $J_i:J_j=J_i+(x^{a_i-a_j})$.

 We begin by establishing a necessary condition for condition (4) (Lemma \ref{necc}), and a sufficient condition for condition (1) (Lemma \ref{suff}).
\begin{lemma}\label{necc}
Let $I=J_1\cap J_2\cap J_3$ be a monomial type three ideal. If $x \in R$ is a monomial such that $xJ_i \inc I:(I:J_i)$ for all $i\in \{1, 2, 3\}$, then $x \in  (I:J_1)+(I:J_2)+(I:J_3)$.
\end{lemma}

\begin{proof}
Assume that $x \notin (I:J_1)+(I:J_2)+(I:J_3)$, but 
for all $i=1, 2, 3$ $xJ_i \inc I:(I:J_i)$. 

By lemma \ref{genlcase} we have
$$
x \in J_1:J_2+J_1:J_3+I:J_1
$$
$$
x\in J_2:J_1+J_2:J_3+I:J_2
$$
$$
x\in J_3:J_1+J_3:J_2+I:J_3
$$
Due to the fact that $x$ and all the ideals involved are monomial, we must in fact have:
$$
x \in J_1:J_2+J_1:J_3
$$
$$
x\in J_2:J_1+J_2:J_3
$$
$$
x\in J_3:J_1+J_3:J_2
$$
Since $x\notin I:J_1$, we may assume either $x \notin J_2:J_1$ or $x \notin J_3:J_1$. Assume for instance $x \notin J_2:J_1$. The second equation implies $x\in J_2:J_3$, and it follows that $x \notin J_1:J_3$ (otherwise $x$ would be in $I:J_3$). 
Thus the first equation implies $x\in J_1:J_2$, so $x \notin J_3:J_2$ (otherwise $x\in I:J_2$), and thus $x\in J_3:J_1$. 

Combining these results, we have
\begin{equation}\label{cond}
x\in (J_2:J_3) \cap ( J_3:J_1)\cap (J_1:J_2).
\end{equation} 
Also note that we cannot have $x \in J_1+J_2+J_3$, since for instance $x\in J_1$ would imply $x \in J_1:J_3$.

Similarly, the assumption $x \notin J_3:J_1$ leads to 
$x\in (J_3:J_2)\cap (J_2:J_1)\cap (J_1:J_3)$.
and $x \notin J_1+J_2+J_3$.

Each of these situations is impossible: for instance condition (\ref{cond}) implies that $x=x^{a_1-a_2}=x^{a_2-a_3}=x^{a_3-a_1}$; if $x_k$ is a variable which appears in $x$ with exponent $c_k>0$, we must have $c_k=a_{1k}-a_{2k}=a_{2k}-a_{3k}=a_{3k}-a_{1k}$, which is clearly impossible.
\end{proof}

On a related note, we have the following general fact:
\begin{prop}
Let $I=J_1\cap J_2\cap \ldots \cap J_n$ be an $m$-primary ideal, with $J_1, \ldots, J_n$ $m$-primary Gorenstein ideals. 
Then
$$\omega ^*(\omega )\inc \frac{I:J_1+\ldots +I:J_n}{I}.$$
\end{prop}
\begin{proof}
Pick $J\inc I$ a Gorenstein $\m-$ primary ideal. According to (~\ref{canonical}),  
$$
\omega \cong \frac{J:J_1+\ldots +J:J_n}{J}.
$$
Since $J:J_i$ is annihilated by $J_i$ in $\omega$, its image under any $f\in \omega ^*$ is also annihilated by $J_i$ in $S/I$, and thus it is contained in $(I:J_i)/I$.
\end{proof}

\begin{lemma}\label{suff}
Let $I=J_1\cap J_2\cap \ldots \cap J_n$ be a type $n$ ideal, with $J_1, \ldots, J_n$ $\m$-primary Gorenstein ideals. Pick $J\inc I$ an $\m-$ primary Gorenstein ideal, and write $J_i=J:f_i$ for some $f_i \in R$. 
If $u_1, \ldots, u_n\in S$ are such that $u_if_j=u_jf_i$ (mod $J$) for all $i, j=1, \ldots, n$, then there is an $R$-linear function $\phi:\omega \rightarrow R=S/I$ defined by $\phi(f_i)=u_i$.
\end{lemma}
\begin{proof}
If $\alpha _1f_1+\ldots \alpha _nf_n \in J$ is a relation on $f_1, \ldots, f_n$ as elements in $\omega$, we must show that $\alpha _1u_1+\ldots +\alpha _nu_n \in I$ is a relation on the images in $S/I$. But this is clear, since
$f_i(\alpha _1u_1+\ldots +\alpha _nu_n )=u_i(\alpha _1f_1+\ldots \alpha _nf_n )=0$ (mod $J$).
\end{proof}
\begin{theorem}
Let $S=k[x_1, \ldots, x_n]$, and let  $I=J_1\cap J_2 \cap J_3$ be a type three $\m$-primary monomial ideal, with $J_1=(x_1^{a_{11}}, \ldots, x_n^{a_{1n}}),
J_2=(x_1^{a_{21}}, \ldots, x_n^{a_{2n}}), J_3=(x_1^{a_{31}}, \ldots, x_n^{a_{3n}})$.
Then the three conditions in Proposition ~\ref{almgor} are equivalent for $R=S/I$, and they hold if and only if $R$ is isomorphic to the ring obtained in one of the following situations:

1. $I$ is a Teter ideal, i.e. $I=(x_1^{a_1}, \ldots, x_n^{a_n}, x_1^{a_1-1}\cdots x_n^{a_n-1})$.

2.  
$$
J_1=(x_1^a, x_2^b, x_3^{c_3}, \ldots, x_n^{c_n})
$$
$$
J_2=(x_1^{a+1}, x_2, x_3, \ldots, x_n)
$$
$$
J_3=(x_1, x_2^{b+1}, x_3, \ldots, x_n).
$$
with $n\ge 3$ and $c_3, \ldots, c_n>1$.

3. 
$$
J_1=(x_1^a, x_2^{b+1}, x_3, \ldots, x_n),
$$
$$
J_2=(x_1^{a+1}, x_2^b, x_3, \ldots, x_n),
$$
$$
J_3=(x_1, x_2, x_3^{c_1}, \ldots, x_n^{c_n}).
$$

4.
$$
J_1=(x_1^a, ,x_2, \ldots, x_s, x_{s+1}^{b_{s+1}}, \ldots, x_n^{b_n}) 
$$
$$
J_2=(x_1^{a+1}, x_2, \ldots, x_n)
$$
$$
J_3=(x_1, x_2^{c_2}, \ldots, x_s^{c_s}, \ldots, x_{s+1}, \ldots, x_n).
$$

5. 
$$
J_1=(x_1, \ldots, x_s, x_{s+1}, \ldots, x_t, x_{t+1}^{c_{t+1}}, \ldots, x_n^{c_n}),
$$
$$
J_2=(x_1, \ldots, x_s, x_{s+1}^{b_{s+1}}, \ldots, x_t^{b_t}, x_{t+1}, \ldots, x_n),
$$
$$
J_3=(x_1^{a_1}, \ldots, x_s^{a_s}, x_{s+1}, \ldots, x_t, x_{t+1}, \ldots, x_n).
$$
\end{theorem}

\begin{proof} First let us check that each of the situations listed in 1.--5. above yields an ideal $I=J_1\cap J_2\cap J_3$ for which condition (1) in Prop. \ref{almgor} holds.

 Consider the five cases in the statement of the theorem.
 Case 1. is covered by Teter's result.

 Consider case 2. Take $J=(x_1^{a+1}, x_2^{b+1}, x_3^{c_3}, \ldots, x_n^{c_n})$, $f_1=x_1x_2$, $f_2=x_2^bx_3^{c_3-1}\cdots x_n^{c_n-1}$, $f_3=x_1^ax_3^{c_3-1}\cdots x_n^{c_n-1}$, so that $J_i=J:f_i$.

For each $i\ge 3$ define $\phi _i:\omega_R\rightarrow R$ by $\phi _i(f_1)=x_i$, $\phi _i(f_2)=x_2^b$, $\phi _i(f_3)=x_1^a$. This is well-defined by Lemma~\ref{suff}, since we have:

 $x_if_2=x_2^bf_1=0$ (mod $J$), $x_if_3=x_1^af_1=0$ (modulo J), and $x_1^af_2=x_2^bf_3$. Thus, $x_3, \ldots, x_n \in \omega^*(\omega)$.

Also define $\phi_i:\omega \rightarrow R$ by $\phi _1(f_1)=x_1,\,  \phi_i(f_2)=x_2^{b-1}x_3^{c_3}-1\cdots x_n^{c_n-1}, \,$ $ \phi_i(f_3)=0$ and $\phi _2:\omega _R\rightarrow R$ by 
$\phi_2(f_1)=x_2, , \phi_2(f_2)=0,$ $\phi_2(f_3)=x_1^{a-1}x_3^{c_3-1}\cdots c_n^{c_n-1}$. It is not hard to check that $x_i\phi(f_j)=x_j\phi(f_i)$ (mod $J$) for all $i, j=1, 2, 3$, where $\phi =\phi _1, \phi_2$, and thus $\phi_1, \phi_2$ are well-defined, and $x_1, x_2 \in \omega^*(\omega)$.

 Consider case 3. Take $J=(x_1^{a+1}, x_2^{b+1}, x_3^{c_3}\ldots, x_n^{c_n})$, $f_1=x_1x_3^{c_3-1}\cdots x_n^{c_n-1}$, $f_2=x_2x_3^{c_3-1}\cdots x_n^{c_n-1}$, $f_3=x_1^ax_2^b$, so that $J_i=J:f_i$, and $\omega =(J, f_1, f_2, f_3)/J$.
For each $i\ge 3$ define $\phi_i:\omega \rightarrow R$ by
$\phi_i(f_1)=x_1, \phi_i(f_2)=x_2, \phi_i(f_3)=x_i$. We have $x_1f_2=x_2f_1$, $x_if_1=x_1f_3=0$ (mod $J$), $x_if_2=x_2f_3=0$ (mod $J$), and thus $\phi _i$ is well-defined by Lemma ~\ref{suff}.

 Consider case 4. Take $J=(x_1^{a+1}, x_2^{c_2}, \ldots, x_s^{c_s}, x_{s+1}^{b_{s+1}}, \ldots, x_n^{b_n})$, $f_1=x_1x_2^{c_2-1}\cdots x_s^{c_s-1}$, $f_2=x_2^{c_2-1}\cdots x_n^{b_n-1}$, $f_3=x_1^ax_{s+1}^{b_{s+1}-1}\cdots x_n^{b_n-1}$. For each $i=2, \ldots, s$ and $j=s+1, \ldots, n$ define $\phi_{ij}:\omega \rightarrow R$ by $\phi_{ij}(f_1)=x_j, \phi_{ij}(f_2)=x_1^a$, $\phi_{ij}(f_3)=x_i$. 
We have $x_1^af_1=x_jf_2=0$ (mod $J$), $x_1^af_3=x_if_2=0$ (mod $J$), and 
$x_if_1=x_jf_3=0$ (mod $J$), thus $\phi_{ij}$ is well-defined, and $x_2, \ldots, x_n \in \omega ^*(\omega)$.

 Also define $\phi:\omega \rightarrow R$ with $\phi(f_1)=x_1$, $\phi(f_2)=x_{s+1}^{b_{s+1}-1}\cdots x_n^{b_n-1}$, $\phi(f_3)=0$. We have $x_{s+1}^{b_{s+1}-1}\cdots x_n^{b_n-1}f_1=x_1f_2$, $x_{s+1}^{b_{s+1}-1}\cdots x_n^{b_n-1}f_3=x_1f_3=0$ (mod $J$), thus $\phi $ is well-defined, and $x_1\in \omega ^*(\omega)$.

 Consider case 5. Take $J=(x_1^{a_1}, \ldots, x_n^{c_n})$,
$f_1=x_1^{a_1-1}\cdots x_t^{b_t-1},
f_2=x_1^{a_1}\cdots x_s^{a_s-1}x_{t+1}^{c_{t+1}-1}\cdots x_n^{c_n-1}
f_3=x_{s+1}^{b_{s+1}-1}\cdots x_n^{c_n-1}$.

For each $i=1, \ldots, s$, $j=s+1, \ldots, t$, $k=t+1, \ldots, n$, we can define a map
$\phi _{ijk}:\omega \rightarrow R$ by
$\phi _{ijk}(f_1)=x_k, \phi _{ijk}(f_2)=x_j, \phi _{ijk}(f_3)=x_i$. We have $x_kf_2=x_kf_3=x_jf_1=x_jf_3=x_if_1x_if_2=0$ (mod $J$), and therefore $\phi_{ijk}$ is well-defined, and $x_1, \ldots, x_n \in \omega^*(\omega)$.

 Now we show  condition (4) in Prop. ~\ref{almgor} implies one of the cases 1.--5. listed in the statement of the theorem. 

Note that whenever we have a variable $x_k \in J_i:J_j$, this implies either $x_k \in J_i$, or else $x_k=x^{a_i-a_j}$. The latter is equivalent to 
$$
\left\{ \begin{array}{c}
a_{ik}=a_{jk}+1 \\
a_{il}\le a_{jl}\, \forall \, l \ne k\\
\end{array}\right.
$$

According to Lemma \ref{necc}, every variable $x_k$ is in one of the ideals $I:J_1$, $I:J_2$, $I:J_3$. 
Without loss of generality we may assume $x_1 \in I:J_1$.
This implies that one of the following holds:

$x_1 \in J_2 \cap J_3$, or
$x_1=x^{a_2-a_1} \in J_3$, or $x_1=x^{a_3-a_1}\in J_2$, or $x_1=x^{a_3-a_1}=x^{a_2-a_1}$.

We claim that $x_1=x^{a_3-a_1}=x^{a_2-a_1}$ implies that $I$ is a Teter ideal.
We have $a_{31}=a_{21}=a_{11}+1$, and $a_{3k}, a_{2k}\le a_{1k}$ for all $k\ne 1$, and therefore we may assume $x_k \notin J_1$ for all $k\ne 1$.
On the other hand, $x_k \in I:J_i$ for some $i$. If $i=1$, it follows that $x_k \in J_2\cap J_3$, because otherwise we would have $x_k=x_1$. 

Since $I$ has type three, there must be at least two variables, say $x_2, x_3$, not in $J_2\cap J_3$ (otherwise we would have a containment between $J_2, J_3$). Without loss of generality, assume $x_2=x^{a_1-a_2}$ (in $I:J_2$) and $x_3=x^{a_1-a_3}$ (in $I:J_3$). It follows that $a_{1k}=a_{2k}$ for all $k\ne 1, 2$, and $a_{1k}=a_{3k}$ for all $k\ne 1, 3$. In particular, all the variables $x_k$ with $k>3$ have $a_{1k}=a_{2k}=a_{3k}$, and therefore they must be in $J_1\cap J_2\cap J_3$ and may be omitted. Therefore we are in case 1.

 Assume that $I$ is not a Teter ideal, thus $x_k \in I:J_i$ implies that $x_k \in J_j+J_l$, where $\{i, j, k\}=\{1, 2, 3\}$.

Assume that we have two distinct variables, say $x_1, x_2 \in I:J_1$, which belong to at most one of $J_1, J_2, J_3$. We must have $x_1=x^{a_2-a_1}\in J_3$ and $x_2 =x^{a_3-a_1} \in J_2$, or vice versa.
For $k>2$ we have $a_{2k}, a_{3k}\le a_{1k}$. If $x_k \in I:J_1$, then $x_k \in J_2\cap J_3$. If $x_k \in I:J_2$, then $x_k=x^{a_1-a_2}$ (otherwise, $x_k \in J_1 \Rightarrow x_k \in J_1\cap J_2\cap J_3$, and we may omit the variable $x_k$ from the presentation of $S/I$). But this implies that $a_{1l}\le a_{2l}$ for all $l\ne k$, contradicting $x_2\in J_2 \, \backslash\,  J_1$. We are in case 2.

 From now on, assume that for each $i\in \{1, 2, 3\}$ there is at most one variable in
$I:J_i$ which is not in the intersection of two of the ideals $J_1, J_2, J_3$.

 Consider the case when each of $I:J_1, I:J_2$ contains a variable which belongs to at most one of $J_1, J_2, J_3$. Without loss of generality, assume that these variables are $x_1 \in I:J_1$, $x_2\in I:J_2$.
 
 If $x_1=x^{a_3-a_1} \in J_2$ and $x_2=x^{a_1-a_2}\in J_3$, we have $a_{1k}\le a_{2k}$ for all $k\ne 2$, and on the other hand $a_{21}=1$, thus $a_{11}=1$, contradicting $x_1 \notin J_1$.

 If $x_1=x^{a_3-a_1} \in J_2, x_2=x^{a_3-a_2} \in J_1$, we have $a_{3k}\le a_{1k}$ for all $k\ne 1$, and on the other hand $a_{12}=1$, thus $a_{32}=1$, contradicting $x_2\notin J_3$.

 If $x_1=x^{a_2-a_1}\in J_3$ and $x_2=x^{a_3-a_2}\in J_1$, we have $a_{2k}\le a_{1k}$ for all $k\ne 1$, and on the other hand $a_{12}=1$, thus $a_{22}=1$, contradicting $x_2\notin J_2$.

 If $x_1=x^{a_2-a_1}\in J_3$, and $x_2=x^{a_1-a_2}\in J_3$, it follows that $a_{1k}=a_{2k}$ for all $k>2$. Thus, $x_k \in J_1:J_2$ or $x_k \in J_2:J_1$ implies $x_k \in J_1\cap J_2$. If $x_k \in I:J_3$, it follows that $x_1x_k\in J_2\cap J_3$, and since $x_1\notin (J_2+J_3)$ we must have $x_k \in J_2\cap J_3$, and thus $x_k \in J_1\cap J_2\cap J_3$, and we may omit these variables. Therefore we are in case 3.

 Now assume that there is only one variable that belongs to exactly one of $J_1, J_2, J_3$, say $x_1 =x^{a_2-a_1}\in J_3$, all other variables belong to exactly two of $J_1, J_2, J_3$. Since $a_{2k}\le a_{1k}$ for all $k\ne 1$, we have either $x_k \in J_2\cap J_1$, or $x_k \in J_2\cap J_3$.
Thus we have case 3.

 Finally, if every variable belongs to exactly two of $J_1, J_2, J_3$, we are in case 5.
\end{proof}


\begin{thebibliography}{99}

\bibitem[AB]{AB}
M.~Auslander and R.-O. Buchweitz,
{\em The homological theory of maximal {C}ohen-{M}acaulay approximations\/},
M\'em. Soc. Math. France (N.S.) (1989),
no.~38, Colloque en l'honneur de Pierre Samuel (Orsay, 1987), 5--37.

\bibitem[D]{D}
S. Ding
{\em A note on the index of Cohen-Macaulay local rings\/},
Communications Alg. {\bf 21} (1993), 53--71.

\bibitem[Te]{Te}
W. Teter
{\em Rings which are a factor of a Gorenstein ring by its socle\/},
Inventiones Math. {\bf 23} (1974), 153--162.

\end{thebibliography}
\end{document}